\magnification=\magstep1    
\hsize=5.9truein
\voffset=.5truein
\input amssym.def
\input amssym.tex
\def\newline{\hfill\break}
\def\scong{{\scriptstyle\|}\lower.2ex\hbox{$\wr$}}
\def\Z{{\Bbb Z}}

\def\Hom{\mathop{\rm Hom}\nolimits}
\def\Ind{\mathop{\rm Ind}\nolimits}

\def\Gal{\mathop{\rm Gal}\nolimits}

\def\rtimes{\mathop{\times\!\!{\raise.2ex\hbox{$\scriptscriptstyle|$}}}
	\nolimits} 
\def\proof{\noindent{\it Proof.}\quad}
\def\blackbox{\hbox{\vrule width6pt height7pt depth1pt}} 
\outer\def\Demo #1. #2\par{\medbreak\noindent {\it#1.\enspace}
	{\rm#2}\par\ifdim\lastskip<\medskipamount\removelastskip
	\penalty55\medskip\fi}
\def\qed{~\hfill\blackbox\medskip}
\overfullrule=0pt

\def\hangbox to #1 #2{\vskip1pt\hangindent #1\noindent \hbox to #1{#2}$\!\!$}
 
\pageno=0

\footline{\ifnum\pageno=0\hfill\else\hss\tenrm\folio\hss\fi}
\topinsert\vskip1.8truecm\endinsert
\centerline{\bf Finite $u$ Invariant and Bounds on Cohomology Symbol Lengths}
\vskip6pt
$${\vbox{\halign{\hfil\hbox{#}\hfil\qquad&\hfil\hbox{#}\hfil\cr$$
$$David J. Saltman\cr
Center for Communications Research\cr
805 Bunn Drive\cr
Princeton, NJ 08540\cr}}}$$
\vskip16pt
\vskip3pt
{\narrower\smallskip\noindent
\hskip 1.8in {\bf{Abstract}}
\bigskip 
In this note we answer a question of Parimala's, showing that fields with finite 
$u$ invariant have bounds on the symbol lengths in their $\mu_2$ cohomology 
in all degrees.

\bigskip

\noindent AMS Subject Classification:  12G05, 11E04
\medskip

\noindent Keywords: u invariant, symbol length\smallskip}
\footnote{}{*The author is grateful for the hospitality of AIM at their workshop 
titled ``Deformation Theory, patching, quadratic forms, and the Brauer group'', 
while this work began.}
\vfill\eject

\leftline{\bf{Introduction}} 
\medskip
At the AIM Workshop on Period/Index problems in January 2011,\break Prof. 
Parimala asked whether fields of finite $u$ invariant necessarily 
had bounded symbol length in their $\mu_2$ cohomology. Parimala 
presented a proof that this was true in degrees one through three, 
using generic splitting constructions. In the subsequent breakout 
session on this problem, the first ideas of a proof were presented 
by myself, and I greatly benefited from the constructive comments by 
Prof. Parimala and Prof. Merkurjev. The note below is my write up 
of the argument. 

Much of the notation and definitions below are standard in quadratic form 
books (e.g. [Sc]), and we also assume familiarity with 
group and Galois cohomology and the Hochschild Serre spectral sequence 
(e.g. [NSW]). As notation which is perhaps not standard, 
for any field $F$ let $G_F$ be its absolute 
Galois group. Let me also add an extended discussion about Galois extensions which is also less standard ([Sa] p.~253). 
Let $S/R$ be an $H$ Galois extension of commutative rings, where 
$H$ is a subgroup of the finite group $G$. Then $\Hom_{H}(\Z[G], S)$ 
can be given the structure of an $R$ algebra by pointwise operations 
in $S$, and has the coinduced $G$ actions. Together this defines 
$\Ind_H^G(S/R)$ which is a $G$ Galois extension of $R$. Note that 
$\Ind_H^G(S/R)$ is of the form $S \oplus \ldots \oplus S$ as an $R$ 
algebra. We will most often use this in the case $R = F$ is a field. 
Then recall that any $G$ Galois $L/F$ ($L$ is not necessarily a field) 
has the form $\Ind_H^G(K/F)$ where $K/F$ is an $H$ Galois extension 
of fields. Finally, if $S$ is a ring on which the group $B$ acts 
we denote by $S^B$ the subring of $B$ fixed elements. 

Underlying this work is the remarkable result of [OVV], based on the
underlyng work of Voevodsky, that shows the maps $I_K^i/I_K^{i+1} \to
H^i(K,\mu_2)$ are isomorphisms for all $i$. Since the paper [OVV] assumes
all fields have characteristic 0 we also make this assumption, but
really this work applies to any fields (characteristic not 2) where
Milnor's Conjecture holds.  We set $\mu_2$ to be the group $\{1,-1\}$
of $2$ roots of 1.  Given a field $F$, then $H^i(F,\mu_2) = H^i(\bar
G,\mu_2)$ the Galois cohomology group where $\bar G$ is the absolute
Galois group of $F$.  If $a \in F^*$ then we abuse notation and write
$a \in H^1(F,\mu_2)$ to be the character $\Hom(\bar G,\mu_2) =
H^1(F,\mu_2)$ determined by the field extension $F(a^{1/2})/F$. The
cup product $a_1 \cup \ldots \cup a_i \in H^i(F,\mu_2)$ is called a
{\bf symbol} and the symbol length of an element $\alpha \in
H^i(F,\mu_2)$ is the least $i$ such that $\alpha$ is a sum of $i$
symbols.  The canonical map $I_K^i/I_K^{i+1} \to H^i(K,\mu_2)$ of
Milnor's Conjecture is determined by mapping the Pfister form
(e.g., [Sc] p.~72) $\langle\langle a_1,\ldots,a_i\rangle\rangle \to -a_1
\cup \ldots \cup -a_i$ and this was shown in [OVV] to be an
isomorphism. In particular, every element of $H^i(F,\mu_2)$ is a sum
of symbols and so has a symbol length.  Recall that the $u$ invariant
of a field is the integer or $\infty$, $u(F)$, such that any quadratic form
over $F$ of rank bigger than $u(F)$ is isotropic. 

If $L = \Ind_H^G(K/F)$, $K$ is a field, and $\beta 
\in H^i(G,\mu_2)$ then $\beta$ has a natural image in Galois cohomology 
by first forming the restriction $\beta_H \in H^i(H,\mu_2)$ and then 
taking the image of $\beta_H$ in $H^i(F,\mu_2) = H^i(G_F,\mu_2)$ via 
inflation. Note that this map commutes with inflation. That is, 
If $G' \to G$ is defined by a tower $F \subset L \subset L'$ of Galois 
extensions, and $L' = \Ind_{H'}^{G'}(K'/F)$ then we can choose 
$H'$, which is only defined up to conjugacy, such that 
we have the diagram: 
$$\matrix
{H'&\subset&G'\cr
\downarrow&&\downarrow\cr
H&\subset&G\cr}$$
and the restriction of $\beta$ to $H$ inflated to $H'$ is the same as 
the inflation of $\beta$ to $G'$ restricted to $H'$. 
\bigskip
\leftline{{\bf{The result.}}}

We say that a field $F$ is {\bf u-bounded} if there is an integer function 
$N(n)$ such that $u(L) \leq N(n)$ for all extensions $L/K$ of degree dividing $n$. 
Note that this is equivalent (e.g., [Sc] p.~104) to just saying that 
$u(F) < \infty$ but we phrase it in this way to emphasize that we are 
considering properties closed under finite extension. 
We say that $F$ has bounded symbol length in degree $d$ 
if there is an integer $M_d(n)$ such that every element in $H^d(L,\mu_2)$ is a 
sum of $M_d(n)$ symbols for every $L/F$ finite of degree dividing 
$n$. The point of this note is to show: 

\proclaim Theorem 1. Every field with finite $u$ invariant has bounded symbol length 
in degree $i$ for all $i$. 

We prove this by induction and we note that every field has bounded symbol 
length of degree 1. Also, since the premise of this result is preserved 
by finite extensions, we may assume we have shown for such $F$ 
that $M_j(n)$ exists for all $j < d$ and $n$, and show $M_d(1)$ 
exists. That is, we show that every element of $H^d(F,\mu_2)$ 
has bounded symbol length.  

The idea of the argument is the following. Since $u(F)$ is finite, we can 
write down a generic quadratic form which specializes to all anisotropic quadratic 
forms. We would like to modify this generic form so that it is generic and lies 
in the Witt Ring fundamental ideal $I^i$ and specializes to all Witt classes of forms in 
that ideal. That form maps to $H^i(F,\mu_2)$ and there it is a sum of some 
number of symbols. By specializing, all elements of $H^i(F,\mu_2)$ 
are the sum of that many or fewer symbols, and we would be done. 

Unfortunately, it does not seem possible to construct such a generic form. 
The difficulty can be illustrated as follows. Suppose 
$\alpha = \sum_{j=1}^m a_{1,j} \cup \ldots \cup a_{i,j}$ is an element of 
$H^i(F,\mu_2)$. Let $L = F(a_{k,l}^{1/2}|\,\, $all$\,\, k,l)$ and let $A$ be the 
Galois group of $L/F$ and $\bar G$ the absolute Galois group of $F$. 
Then our form for $\alpha$ defines an element of $H^i(A,\mu_2)$. If this element maps to 
$0$ in $H^i(\bar G,\mu_2) = H^i(F,\mu_2)$, then there must be an finite extension 
$L' \supset L$ with $L'/F$ Galois with group $B$ and so an induced surjection $B \to A$ 
such that $\alpha$ maps to 0 in $H^i(B,\mu_2)$. However, we see no way to, 
in general, bound the size of $B$ and if the size of $B$ is not bounded 
there can be no generic way to force a cohomology class of $\alpha$'s 
form to be 0 because the $B$ that works for such a generic construction 
would then give you a bound. However, we can construct a generic zeroing of 
$\alpha$ for a fixed $B$. We will call this generic with the limitation 
$B$ (formal definition to follow). Before we explore this, we mention the 
following way of thinking about writing a cohomology element as a sum of symbols. 

\proclaim Lemma 2. An element of $\alpha \in H^i(F,\mu_2)$ is a sum of symbols 
if and only if $\alpha$ is the image of $H^i(A,\mu_2)$ where $A = \Gal(L/F)$ 
is an elementary abelian $2$ group. Moreover, that all $\alpha$ have bounded 
symbol length is equivalent to bounding the size of such $A$. 

\proof The basic equivalence is immediate from (e.g., [E] p.~33), which 
says the well-known fact that any element of such an $H^i(A,\mu_2)$ is the 
sum of $i$ degree monomials of elements in $H^1(A,\mu_2) = \Hom(A,\mu_2)$. 
Moreover, if $\alpha = \sum_{j=1}^M (a_{j,1}\cup\ldots\cup{a_{j,i}})$ then $A$ 
can be taken of order dividing $2^{im}$. Conversely, for $A$ of order $2^M$, 
$H^i(A,\mu_2)$ has a basis consisting of ${M}\choose{i}$ elements.~\qed 

In what follows $A$ will always be an elementary abelian 2 group. 

Let $C$ be our ground field so all rings and fields will be $C$ algebras. 
We suppose $\alpha \in H^i(F,\mu_2)$ is the image of $\beta \in H^i(A,\mu_2)$ 
where $A = \Gal(L/F)$ is as above. Further, suppose $S/R$ is a 
$A$ Galois extension with $q(R) = F$, $R$ affine over $C$, and $L = S \otimes_{R} F$ as 
$A$ Galois extensions. We call $S/R$, $A$ and $\beta$ a {\bf presentation} 
of $\alpha$. Since we are usually not interested in the specific rings 
$S/R$ in a presentation, we define $S'/R'$, $A$, $\beta'$ {\bf equivalent} 
to $S/R$, $A$, $\beta$ if and only if $A = A'$, $\beta = \beta'$, 
and there are nonzero $r \in R$ and $r'\in R'$ such that 
$R(1/r) = R'(1/r')$ and $S(1/r)/R(1/r)$ and $S'(1/r')/R'(1/r')$ 
are isomorphic as $A$ Galois extensions of $R(1/r)$. Obviously 
equivalent presentations have the same induced cohomology element. 
In discussing presentations up to equivalence, the ring extension 
$S/R$ can often be surpressed and we can just say the presentation 
$\beta$, $L/F$ where $L/F$ is $A$ Galois. 

Presentations will be important to us because they allow specializations 
of cohomology classes as follows. 
In fact, we will be defining specializations of presentations. 
Let $\beta$, $S/R$, and $A$ be a presentation of $\alpha$. 
Suppose $\phi: R \to R_1 \subset F_1$ for a ring and field 
$R_1, F_1 \supset C$. 
If we set $S_1 = S \otimes_{\phi} R_1$ then $\beta$, $S_1/R_1$, and $A$ 
is the specialization with respect to $\phi$. 
If $L_1 = S \otimes_{\phi} F_1$ then this specialized presentation 
defines an $\alpha_1 \in H^i(F_1,\mu_2)$ which we can call a 
specialization of $\alpha$. 

Note that we have defined the notion of presentation without assuming 
$L = S \otimes_R F$ is a field. In fact, suppose $L = \Ind_{A_1}^A(L_1/F)$. 
Then there is some $0 \not= r \in R$ such that 
$S(1/r) = \Ind_{A_1}^A(S_1/R(1/r))$ and we can define (up to equivalence) 
$\beta_1$, $S_1/R(1/r)$, $A_1$ to be a restriction of the original 
presentation, and this restriction presents the same cohomology element. 

Next we must talk about presentations that represent 0 
and their so called limitations. 
Let $\beta \in H^i(A,\mu_2)$ be as above and suppose $L/F$ is an $A$ Galois extension. 
We say $\beta$ and $L/F$ {\bf presents} 0 if the image of $\beta$ in $H^i(F,\mu_2)$ 
is $0$. We say that $\beta$ and $L/F$ present $0$ with {\bf limitation} 
$B \to A$ (sometimes we write only $B$) if and only if there is a Galois 
extension $L'/F$ containing $L/F$ and inducing $B \to A$ such that $\beta$ 
maps to $0$ in $H^i(B,\mu_2)$. If $L = \Ind_{A_1}^A(L_1/F)$ and 
$\beta$ restricts to $\beta_1$, and if $\beta$ has limitation $B$ 
then $\beta_1$ also presents 0 and has limitation $B_1$ where 
$B_1 \subset B$ is the inverse image of $A_1$.  

If $\beta$ and $L/F$ presents $0$ and $L$ is a field then 
there is a $L'/F$ Galois containing $L/F$ and an associated surjection of 
Galois groups $B \to A$ such that $\beta$ presents 0 with limitation $B$. 
More generally, if $L = \Ind_{A_1}^A(L_1/F)$ with $L_1$ a field 
there is a Galois $L_1'/F$ containing $L_1/F$ and associated 
surjection $B_1 \to A_1$ 
such that $\beta$ maps to 0 in $H^i(B_1,\mu_2)$. 
Since $A$ is an elementary abelian 2 group $A_1$ is always a direct 
summand of $A$. If $\rho: A \to A_1$ is any homomorphism which 
is the identity on $A_1$, we can set $B$ to be the pullback 
of $B_1 \to A_1$ and $\rho$. Note that $A_1 \subset A$ induces an injection 
$B_1 \subset B$ and if $L = \Ind_{B_1}^B(L_1')$ then 
$\beta$ presents 0 with limitation $B$. Note also that $B_1$ is a direct 
summand of $B$. 
Of course, there is a restriction 
$\beta_1$ defined by $A_1$ and $L_1$ and of course $\beta_1$ presents 0 with limitation $B_1$. Since $B$ determines $B_1$, and $B_1$ is a direct 
summand, we will sometimes say $\beta_1$ has limitation $B$ .  

Note that since the inflation of $\beta$ to $H^i(B,\mu_2)$ 
is independent of any fields, for $\beta$, $L/F$ to present 0 with 
limitation $B$ really means that $\beta$ maps to 0 in $H^i(B,\mu_2)$ and 
$L/F$ extends 
to a $B$ Galois extension $L'/F$ inducing $B \to A$.  
In addition, if $L$ is a field but $L' = \Ind_{B_1}^B(L_1'/F)$ then 
$B_1 \to A$ is surjective (since $L$ is a field) and the image of $\beta$ 
is also 0 in $H^i(B_1,\mu_2)$ meaning that $\beta$ also has limitation $B_1$. 

Frequently when we specialize as above there will be many ways to do 
it and this is important. Given a presentation $S/R$, and $\beta \in H^i(A,\mu_2)$, we 
say $\beta$ {\bf densely specializes} to a presentation 
$\beta_1$ of $\alpha_1 \in H^i(F_1,\mu_2)$ if the following 
holds. For any $0 \not= r \in R$, there is a $\phi: R \to F_1$ such that 
$\phi(r) \not= 0$ and $\phi$ causes $\beta$ to specialize to $\beta_1$ 
inducing the same presentation $\beta_1$. 
If $R$ and $R'$ are affine $C$ algebras with $q(R) = F = q(R')$, then 
$R(1/r) = R'(1/r')$ for some $0 \not= r \in R$ and $0 \not= r' \in R'$ (e.g. [Sw] p.~152). Thus when $\beta$ densely specializes to $\beta_1$ 
this is well defined up to equivalence.

\proclaim Lemma 3. Suppose $\beta \in H^i(A,\mu_2)$ and $L/F$, $A$  are a presentation 
of $0$ which densely specializes to $\beta_1$, $L_1/F_1$, $A$. 
Then if $\beta$ has limitation $B \to A$ so does $\beta_1$ 
and in particular $\beta_1$ presents 0.  

\proof 
Assume $L' \supset L \supset F$ is $B \to A$ Galois and $\beta$ maps to 0 
in $H^i(B,\mu_2)$. If $S/R$ is $A$ Galois and $q(R) = F$, there is a $0 \not= r \in R$ 
and a $S'/R(1/r)$ which is $B$ Galois, contains $S(1/r)$, with 
$S' \otimes_{R(1/r) } F = L'$. Choose $\phi: R(1/r) \to F_1$ realizing the specialization 
and set $\tilde L_1' = S' \otimes_{\phi} F'$. 
Then $L_1' \supset L_1 \supset F_1$ is $B \to A$ Galois.~\qed 

Let $\alpha \in H^i(F,\mu_2)$ have presentation $\beta$, $A = \Gal(L/F)$. 
Assume $B \to A$ is a surjection of finite groups and $\beta$ maps to 0 
in $H^i(B,\mu_2)$. The following result is routine and we only include the 
proof for ease of the reader.

\proclaim Proposition 4. There is a field extension $F_B \supset F$ 
with the following properties (where we set $L_B = L \otimes_F F_B$). 
\medskip
a) The extension $\beta$, $A = \Gal(L_B/F_B)$ presents 0 
with limitation $B$. 
\medskip
b) Suppose $\beta$, $A = \Gal(L/F)$ densely specializes to $\beta_1$ and 
$A_1 = \Gal(K_1/F_1)$. If $\beta_1$ has limitation $B$, then the extension 
$\beta$, $L_B/F_B$ densely specializes to $\beta_1$. 

\proof 
All this really means is that we are constructing generically the $B$ Galois 
extension extending $L/F$. To achieve this let $V$ be the faithful $B$ 
module $F[B]$. Let $B$ act on the field of fractions $L(V)$ as follows. 
$B$ acts on $V$ as usual, and $B$ acts on $L$ via $B \to A$. 
Set $F_B = L(V)^B$. It is clear that $L(V)/F_B$ is $B$ Galois. 
If $N$ is the kernel of $B \to A$, then $L(V)^N = L_B = L \otimes_F F_B$ 
so $L(V) \supset L_B \supset F_B$ induces $B \to A$.  

Let $S/R$ be $A$ Galois such that $q(R) = F$. Then we can 
choose $t \in S[V]^B$ with the property that if $S_B = S[V](1/t)$ 
then $S_B/R_B$ is $B$ Galois with $R_B = (S_B)^B$. 
Suppose $0 \not= s \in S[V]^B$. It suffices to show that there is a 
$\phi: R_B \to F_1$ with $\phi(s) \not= 0$. By assumption there is a 
$\phi: R \to F$ specializing $\beta$ to $\beta_1$. Set $L_1 = S \otimes_{\phi} F_1$ 
which is $A$ Galois over $F_1$ and has the form $\Ind_{A_1}^A(K_1)$. 
It follows that $\phi$ extends to an $A$ morphism $\phi: S \to L_1$. Since 
$\beta_1$ has limitation $B$, there is a $B$ Galois extension $L_1'/F \supset L_1/F$ 
inducing $B \to A$. Since $V$ has basis $\{x_g|g \in B\}$ with the obvious action, 
algebraic independence of Galois group elements 
(e.g. [BAI] p. 294) shows that we can define $\phi(x_g) = g(a) \in L_1'$ 
for some $a$ such that $\phi(st) \not= 0$. Then $\phi$ extends to a $B$ morphism 
and restricts to the needed $\phi$ on $R_B$.~\qed 

If $\beta$ and $B \to A$ are as in Proposition 4, we say that $F_B$ 
is the generic splitting field of $\beta$ with limitation $B$. 

Let's outline our argument a bit. We start with a generic quadratic form 
$\gamma = \sum^N_{i=1} a_ix_i^2$ ($N$ is even) meaning that the ground field has the form 
$F_1 = C(a_1,\ldots,a_N)$ and the $a_i$ are a transcendence base. Note that 
for any field $F \supset C$, this specializes to all Witt classes in the fundamental ideal 
$I$ as long as $u(F) \leq N$. We want to write down a generic element in $I^n$ 
with a fixed so called history as follows. 
$F_2/F_1$ is the extension defined by taking the square root of the 
determinant of $\gamma$. The extension $\gamma_2$ of $\gamma$ to 
the Witt ring $W(F_2)$ is in $I_{F_2}^2$ and so defines an element 
$\alpha_2 \in H^2(F_2,\mu_2)$. We take $F_3/F_2$ to be a generic splitting field 
of $\alpha_2$  and so the extension, $\gamma_3 \in W(F_3)$ is in $I_{F_3}^3$. 
So far there has been no limitations. However, if $\alpha_3 \in H^3(F_3,\mu_2)$ 
is the image of $\gamma_3$ then we can write $\gamma_3$ as a sum of Pfister forms 
and thereby write $\alpha_3$ as a sum of symbols. Given that, we can choose 
a presentation $\beta_3$, $A_3 = \Gal(L_3/F_3)$ of $\alpha_3$. For 
any $B_3 \to A_3$ that splits $\beta_3$, we form the generic splitting field 
of $\beta_3$ with limitation $B_3$ and call that $F_4$. We proceed by induction 
until the extension, $\gamma_n \in I_{F_n}^n$ is defined. The choice 
of presentations $\beta_i$ and limitations $B_i$ is the {\bf history} 
of this construction. 

Now given a $u$ bounded field $K$ every element $\alpha' \in H^i(K,\mu_2)$ 
is the image of a quadratic form $\gamma'$ which is in $I_K^i$. 
We show that we can bound the order and hence number of the limitations 
which enforce this property of $\gamma'$'s and hence write $\alpha'$ 
as the specialization of one of finitely many of the generic contructions 
of $\alpha_n$ (as above) as we vary the histories among finitely many choices 
of the $B_i$. This proves the result. 

To make this argument more formal, if $\beta$, $A = \Gal(L/K)$ is a presentation 
$\alpha$ then the {\bf order} of $\beta$ is the order of the group $A$. 
Obviously the order of a presentation cannot increase under specialization. 
We say that a field $K$ is {\bf limitation bounded} in degree $i$ if and only if 
for all $d$, all field extensions $K'/K$ of degree dividing $d$, and all 
presentations of zero $\beta$ over $K'$ of order less than or equal to 
$N$, there is a $L(N,d)$ such that $\beta$ has a limitation $B$ of order 
less than or equal to $L(N,d)$. The above argument is an outline of the proof of: 

\proclaim Theorem 5. Suppose $K$ is limitation bounded in all degrees 
$j < d$ and is also $u$ bounded. Then all finite extensions of $K$ have bounded 
symbol length in degree $d$ and the bound is a function of the degree alone. 

\proof This is perhaps already clear except for the fact we are choosing 
presentations of zero. For simplicity we only treat $K$ itself, the extension 
to the $K'/K$ being clear. We prove by induction that every $\gamma' \in I_K^d$ is the specialization of 
some $\gamma_i \in I_{F_d}^d$ as above with only finitely many choices 
of histories. By induction there are finitely many histories such that $\gamma'$ 
is the specialization some $\gamma_{d-1}$. For this $\gamma_{d-1}$ there is a presentation 
$\beta_{d-1}$ and thus a degree $d-1$ presentation of $\gamma'$ we call $\beta_{d-1}'$. 
Since $\gamma' \in I_K^d$ it follows that $\beta_{d-1}'$ is a presentation of zero 
and so there are only finitely many further choices of limitations $B_d$. 
We are done by Proposition 4.~\qed 

Given Theorem 5, we need to prove these $u$ bounded fields are limitation bounded. 
This is an 
involved argument using the Hochschild--Serre spectral sequence. Note that we feel  
that the limitation bound we obtain is far from optimal. For this reason we will not 
particularly explicit about the bound, as in the definition of ``predictable'' 
below. However, there is a group structure bound in our argument that seems 
interesting and so we will endeavor to prove it and make it explicit. 
In fact, let $G$ be a finite group (for us usually abelian). A {\bf $d$-abelian $G$ group} 
is an extension $1 \to N \to G' \to G \to 1$ such that $N$ contains $G'$ normal 
subgroups $N = N(0) \supset N(1) \supset \ldots \supset N(d) = 0$ with 
$N(i)/N(i+1)$ abelian. Given any $d$ abelian group we will use obvious modifications 
of the $N(i)$ notation above to denote the corresponding tower of groups.  

Fix a field $K$ with an absolute Galois group $\bar G = \Gal(\bar K/K)$. We say $G'$ is a 
$d$-abelian $G$ Galois group over $K$ to mean $G'$ is also an image of $\bar G$ (so 
$G' = \Gal (L(d)/K)$). Set $L(i)  = L(d)^{N(i)}$ and $L = L(0)$. 
Whenever we talk about $d$ abelian $G$ Galois groups we will use the $L(i)$ 
notation, or obvious variants of it, to indicate the associated tower of fields. 

In the course of the proof 
we will alter $G'$ in several ways. In all cases we will want to construct Galois 
groups so we will usually construct these further groups 
via field theory. 

If $L'/L$ is abelian with 
$L'/K$ Galois, then $L'(d)/K = L(d)L'/K$ is Galois with group $G''$ which is still 
a $d$-abelian $G$ group. To see this, set $L(i)' = L(i)L'$ for $i > 0$ and 
$L'(0) = L(0)$. Then $L'(i)/L'(i-1)$ is abelian Galois with group a subgroup of 
$N(i-1)/N(i)$ for $i > 1$, but $L'(1)/L$ is abelian with Galois group mapping onto 
$N_1/N_0$. 
We call this {\bf expanding} the $d$-abelian $G$ group $G'$. 

Another construction we will need is the following. 
Suppose $K'(d')/K$ is Galois with $d'$ abelian Galois group so $K'(d') \supset \ldots \supset 
K'(1) = L(1) \supset L(0) = K'(0) = L \supset K$ the point being here is that 
the beginning of the series of fields for $K'$ coincides with the beginning for $L(d)$
(and $K'(i)/K'(i-1)$ is, of course, abelian Galois).  
Let $d''$ be the maximum of $d$ and $d'$, 
and set $K'(j) = K'(d')$ for $j > d'$ and similarly $L(j) = L(d)$ for $j > d$. 
Then if $L'(i) = L(i)K'(i)$ we have that $G'' = \Gal(L'(d'')/K)$ is a $d''$-abelian 
$G$ group. Moreover, the first $\Gal(L'(1)/L'(0))$ is unchanged but the rest of the 
abelian series is larger. We call this {\bf refining} the group $G'$. 
Note that expanding $G'$ increases $\Hom(N_0,\mu_2)$ and so increases the 
cohomology cup products in $H^q(N_0,\mu_2)$. On the other hand, we will see that by 
refining our $d$ abelian $G$ groups we will introduce more relations among these 
cup products. 

If $L(d)/K$ and $L_1(d')/K$ have $d$ and $d'$ abelian $G$ Galois groups $G'$ and 
$G_1$respectively (with the same $G$ Galois $L/K$) then the amalgamation $L(d)L'(d')/K$ 
has a $d''$ abelian $G$ Galois group $G_1'$ where $d''$ is the maximum of $d$ and $d'$. 
We call $G_1'$ the {\bf amalgamation} of $G'$ and $G_1$. 

Suppose $A' = \Gal(L'/L)$ is abelian and we have a $d$ abelian $A'$ Galois group 
over $L$, with associated field extensions $L'(d) \supset \ldots \supset L'(0) = L'$. 
Let $L/K$ be $G$ Galois as above. 
Then $L'(d)$ is not Galois over $K$, but if $L''(d)$ is the Galois closure of $L'(d)$ 
over $K$, then $L''(d)$ is the amalgamation of all the $G$ conjugates of $L'(d)$ 
and so $\Gal(L''(d)/K)$ is a $d+1$ abelian $G$ Galois group over $K$. 
We call this {\bf extending} the $d$-abelian $A'$ group to a $d+1$ abelian $G$ group. 

As another bit of notation, if we have $d$ and $d_1$ abelian $G$ groups with a diagram 
$$\matrix{
1&\to&N_1&\to&G_1'&\to&G&\to&1\cr
&&\downarrow&&\downarrow&&||\cr
1&\to&N&\to&G'&\to&G&\to&1\cr}$$
where all vertical arrows are surjective then we say the $d_1$ abelian $G$ group 
$G_1'$ is a {\bf cover} of $G'$ and if these maps are induced by field extensions 
we call it a Galois cover. Clearly, expanding, refining, amalgamating and extending are ways 
of constructing Galois covers. 

When we expand or refine or extend a $d$-abelian $G$ group $G'$ we say that the size of the new 
group is {\bf predictably} bounded if the bound is only a function of $G'$,  
the degrees of the cohomology groups involved, and previously proven symbol length 
bounds for field extensions of bounded degree. Note how unspecific this notion is. 
Any function of predictable bounds, or functions of predictable bounds and 
$|G|$ etc., also would be a predictable bound. For example, when 
we expand, or refine, or extend or amalgamate predictably bounded groups we get other ones.  

Next we need some notation to help us navigate through the complexities 
of the Hochschild-Serre spectral sequence. We will employ this spectral sequence 
for sequences $1 \to \bar N \to \bar G \to G \to 1$ and 
$1 \to \bar N/\bar N' \to G' \to G \to 1$ 
where $\bar N/\bar N'$ is finite. Of course the natural map defines a morphism 
from the second spectral sequence to the first. 
Let me define notation in the first case as extension 
to the second is obvious. 

In this spectral sequence, the $E^{p,q}_2$ term is $H^p(G,H^q(\bar N,\mu_2))$. 
The differential of this spectral sequence is $d_r$ so $d_2: H^p(G,H^q(\bar N,\mu_2)) 
\to H^{p+2}(G,H^{q-2+1}(\bar N,\mu_2))$. We wish to treat each $E^{p,q}_r$ as a subquotient 
of $H^p(G,H^q(\bar N,\mu_2))$ and so write $E^{p,q}_r = 
H^p(G,H^q(\bar N,\mu_2))_r^u/H^p(G,H^q(\bar N,\mu_2))_r^l$. Thus 
$H^p(G,H^q(\bar N,\mu_2))_2^u = 
\break H^p(G,H^q(\bar N,\mu_2))$ and 
$H^p(G,H^q(\bar N,\mu_2))_2^l = 0$. Moreover, the differentials $d_r$ can be viewed as 
morphisms 
$$d_r^{p,q}: H^p(G,H^q(\bar N,\mu_2))_r^u \to 
H^{p+r}(G,H^{q-r+1}(\bar N,\mu_s))/H^{p+r}(G,H^{q-r+1}(\bar N,\mu_s))_r^l$$ 
and the kernel of $d_r^{p,q}$ is $H^p(G,H^q(\bar N,\mu_2))_{r+1}^u$ while the image of 
$d_r^{p,q}$ is 
$$H^{p+r}(G,H^{q-r+1}(\bar N,\mu_s))_{r+1}^l/H^{p+r}(G,H^{q-r+1}(\bar N,\mu_s))_r^l.$$ 
Since all $d_r^{d,0}$ are 0, the kernel of $H^d(G,H^0(\bar N,\mu_2)) \to H^d(\bar G,\mu_2)$ 
is the union of all the $H^d(G,\mu_2)_r^l$ for all $r$. 

In the arguments to follow we are going to use the above notions and work through 
the details of the Hochschild--Serre spectral sequence showing we can stay 
``predictably bounded'' all the way. Since the overall argument is by induction, 
we will be able to assume the following at the degree $d$ cohomology step: 

(*) For all $q < d$ the following two facts hold. First, for all finite extensions 
of $L/K$ of degree dividing $n$ there is a symbol length bound for $H^q(L,\mu_2)$ only depending 
on $n$. Second, for any abelian extension $L'/L$ of predictably bounded degree and Galois group $G_1$, and any 
element $\beta \in H^q(G_1,\mu_2)$ which presents 0, 
there is a predictably bounded $q-1$ abelian group $G_1$ Galois group 
over $L$ written $1 \to N_1 \to G_1' \to G_1 \to 1$ such that $\beta$ maps to 
$0$ in $H^q(G_1',\mu_2)$. 

\proclaim Lemma 6. Suppose $1 \to N \to G' \to G \to 1$ is a $t$-abelian $G$ Galois 
group over $K$ and $\gamma \in H^p(G,H^q(N,\mu_2))$ maps to 0 in $H^p(G,H^q(\bar N,\mu_2))$ 
where $q < d$. Assume (*). 
Let $q'$ be the maximum of $t$ and $q$. Then there is a 
predictably bounded $q'$ abelian cover of $G'$ where $\gamma$ maps to 0. 

\proof The element $\gamma$ is represented by a $p$-cocycle 
$c(g_1,\ldots,g_p) \in H^q(N,\mu_2)$ whose image in $H^q(\bar N,\mu_2))$ is the coboundary 
of a $p-1$ cochain $b(g_1,\ldots,g_{p-1})$. Since each $b(g_1,\ldots,g_{p-1})$ 
is a sum of $M([L:K])$ symbols, there is a $1$-abelian $G$ Galois group 
$1 \to \bar N/\bar N_1 \to G'' \to G \to 1$ such that the $b$'s are the image of elements 
of $H^q(\bar N/\bar N_1,\mu_2)$ where $\bar N/\bar N_1$ is a predictably bounded elementary abelian 2 group. Expanding $G'$ by this $G''$, and calling the result $G'$ again, 
we have that the $c$'s and $b$'s are both in $H^p(G,H^q(N,\mu_2))$. 
There are $|G|^p$ relations in $H^q(N,\mu_2)$ that must be satisfied 
in order that the coboundary $\delta(b)$ equals $c$. Now let $G'' = N/N(1)$ be  
the abelian group. By assumption, we can iteratively refine the $t-1$ abelian 
$G''$ Galois group $1 \to N(1) \to N \to G'' \to 1$ to force these relations, and 
result is a $q'-1$ abelian $G''$ Galois group. Extending this to a $q'$ abelian $G$ 
Galois group we get the cover we need. Note that the size of the new group is again 
predictably bounded, though the bound is quite large.~\qed 

In a similar vein is: 

\proclaim Lemma 7. Assume (*). Suppose $p + q < d$ and $\gamma \in H^p(G,H^q(\bar N,\mu_2))$. 
Then there is a $q$-abelian $G$ Galois group $1 \to N \to G' \to G \to 1$ 
of predictably bounded order and an element $\gamma' \in H^p(G,H^q(N,\mu_2))$ 
which inflates to $\gamma$. 

\proof The element $\gamma$ is represented by a $p$ cocycle $c(g_1,\ldots,g_p)$ 
consisting of less than or equal to $|G|^p$ elements of $H^q(\bar N,\mu_2)$. 
Each of these elements can be written as a sum of $M(|G|)$ symbols so that there are at most  
$|G|^pM(|G|)$ symbols and hence $|G|^2M(|G|)q$ elements of $\Hom(\bar N,\mu_2)$ 
are involved in writing all the $c(g_1,\ldots,g_p)$'s.  Said differently, 
there is a $\bar N' \subset \bar N$ of index dividing $2^{|G|^pM(|G|)q}$ such that 
all these $c(g_1,\ldots,g_{d-r})$'s are in the image of $H^q(\bar N/\bar N',\mu_2)$. 
Now $\bar N'$ is not normal in $\bar G$ so we take the intersection of the $|G|$ 
conjugates of $\bar N'$ to define $\bar N_1$ such that $\bar N_1$ is normal 
in $\bar G$ and all the $c(g_1,\ldots,g_{d-r})$'s come from 
$c_1(g_1,\ldots,g_p) \in H^q(\bar N/\bar N_1,\mu_2)$. Note that $N/N_1$ 
has order dividing $2^{|G|^{p+1}M(|G|)q}$.
Now the $c_1$'s do not form a $q$ cocycle neccessarily, but their image in 
$H^q(\bar N,\mu_2)$ is a cocycle. Being a cocycle means that there are less than or 
equal to $|G|^{p+1}$ relations that must be satisfied.  
Since $q-1 < d$, by assumption (*) 
there is a $q-1$ abelian $G'$ group where each of the cocycle relations 
become true after inflation (of predictably bounded size). 
By extending we get an $q$ abelian 
$G$ group for each needed cocycle relation, and we can refine all these together 
to get an $q$ abelian $G$ group $1 \to N \to G' \to G \to 1$ of predictably bounded 
size such that the $c_1$'s inflate to an element 
$\gamma' \in H^p(G,H^q(N,\mu_2))$ which inflates to $\gamma$.~\qed 

The previous results do not suffice, as we need to show that we can 
achieve, element by element, the spectral sequence filtration 
in a predictably bounded cover. This is the next result. 

\proclaim Lemma 8. Suppose $p + q < d$, and $1 \to N \to G' \to G \to 1$ 
is a $t$ abelian $G$ Galois group over $K$ and $\gamma' \in H^p(G,H^q(N,\mu_2))$ 
maps via inflation to an element $\gamma \in H^p(G,H^q(\bar N,\mu_2))^u_s$. 
Let $t'$ be the maximum of $t$ and $q-1$. Assume (*). 
Then there is a $t'$ abelian Galois cover $1 \to N_1 \to G_1 \to G \to 1$ 
of predictably bounded size such that the inflation of $\gamma'$ in 
$H^p(G,H^q(N_1,\mu_2))$ lies in $H^p(G,H^q(N_1,\mu_2))^u_s$. 

\proof 
We prove this by induction on $s$. The statement is vacuous for $s = 2$ 
and by way of illustration the fact that $\gamma'$ lies in $H^p(G,H^q(N_1,\mu_2))^u_3$ 
is equivalent to $d_2(\gamma') = 0 \in H^{p+2}(G,H^{q-1}(N_1,\mu_2))$ 
which (by Lemma 6) we can achieve after refining and thereby creating a predictably bounded $t'$-abelian 
$G$ Galois cover where $t'$ is the maximum 
of $t$ and $q-1$. 

So assume the result for $s-1$. We need to unpack the meaning when we say 
$\gamma \in H^p(G,H^q(\bar N,\mu_2))^u_s$. By definition, this is equivalent to 
$d_{s-1}(\gamma) \in H^{p+s-1}(G,H^{q-s+2}(\bar N,\mu_2))_{s-1}^l$ 
or $d_{s-1}(\gamma) - d_{s-2}(\gamma_1) \in H^{p+s-1}(G,H^{q-s+2}(\bar N,\mu_2))^l_{s-2}$ 
where $\gamma_1 \in H^{p+1}(G,H^{q-1}(\bar N,\mu_2))_{s-2}^u$. 
Proceeding by induction we have elements 
$\gamma_i \in H^{p+i}(G,H^{q-i}(\bar N,\mu_2))^u_{s-i-1}$ 
for bounded $i$ such that 
$d_{s-1}(\gamma) = \sum_i d_{s-i-1}(\gamma_i) \in H^{p+s-1}(G,H^{q-s+2}(\bar N,\mu_2))$. 
By repeated use of Lemma 7 the $\gamma_i$ are the image of $\gamma_i' \in 
H^{p+i}(G,H^{q-i}(N,\mu_2))$ for a $q-1$ abelian $G$ group $G'$. By induction 
we can assume the $\gamma_i' \in H^{p+i}(G,H^{q-i}(N,\mu_2))^u_{s-i-1}$ 
and one further refinement allows us to assert that 
$d_{s-1}(\gamma') = \sum_i d_{s-i-1}(\gamma_i') \in 
H^{p+s-1}(G,H^{q-s+2}(N,\mu_2))^l_{s-1}$.~\qed 

Until now our spectral sequence notation has been unambiguous as to whether 
we are referring to the absolute sequence $1 \to \bar N \to \bar G \to G \to 1$ 
or some finite image of it, but we have to make such a distinction when 
we deal with $H^d(G,\mu_2)$ which appears in all these spectral sequences. 
Thus we will let $H^d(G,\mu_2)^l_r$ refer to the filtration induced 
by the absolute sequence and for the image sequence 
$1 \to \bar N/N' \to G' \to G \to 1$ we will use the notation $H^d_{G'}(G,\mu_2)^l_r$. 

\proclaim Proposition 9. Assume (*). Suppose $\beta \in H^d(G,\mu_2)_{r+1}^l$ for some $r$. 
Then there is a $r-1$ abelian $G$ Galois group $1 \to N \to G' \to G \to 1$ 
of predictably bounded order such that $\beta \in H^d_{G'}(G,\mu_2)_{r+1}^l$. 

\proof There is a 
$\gamma \in H^{d-r}(G,H^{r-1}(\bar N,\mu_2))_r^u$ 
which maps to $\beta$ modulo $H^d(G,\mu_2)_r^l$. Let $\eta = d(\gamma) - \beta 
\in H^d(G,\mu_2)_r^l$. By induction on $r$ there is a predictably bounded 
$r-2$ abelian $G$ Galois group $1 \to N_1 \to G_1 \to G \to 1$ 
such that $\eta$ is the inflation of $\eta' \in H^d_{G_1}(G,\mu_2)_r^l$.
By Lemma 7 there is a predictably bounded 
$r - 1$ abelian $G$ Galois group $1 \to N \to G' \to G \to 1$ such that 
$\gamma$ is the inflation of $\gamma' \in H^{d-r}(G,H^{r-1}(N,\mu_2))$. 
Since $\gamma \in H^{d-r}(G,H^{r-1}(\bar N,\mu_2))_r^u$ we know from 
Lemma 8 that there is $r-1$ abelian $G$ Galois group cover of $G'$, 
which we also call $G'$, such that $\gamma' \in H^{d-r}(G,H^{r-1}(N,\mu_2))_r^u$. 
We can amalgamate $G_1$ and $G'$ to get a predictably bounded $r-1$ abelian $G$ 
Galois group (we also call $G'$) where $d_r(\gamma')$, $\beta$ and $\eta'$ are all 
defined and $\beta = d_r(\gamma') - \eta'$ so $\beta \in 
H^d_{G'}(G,\mu_2)^l_{r+1}$.~\qed 

\proclaim Corollary 10. Suppose $\beta \in H^d(G,\mu_2)$ maps to 0 
in $H^d(\bar G,\mu_2)$. Assume (*). Then there is a predictably bounded $d-1$ abelian $G$ Galois group 
such that $\beta$ maps to 0 in $H^d(G',\mu_2)$. 

\proof In the spectral sequence the last nontrivial derivation with image $H^d(G,\mu_2)$ 
is $d_d: H^0(G,H^{d-1}(\bar N,\mu_2))^u_d \to H^d(G,\mu_2)/H^d(G,\mu_2)^l_d$. 
That is, the ascending tower in $H^d(G,\mu_2)$ stabilizes at $H^d(G,\mu_2)^l_{d+1}$. 
The result follows from Proposition 9.~\qed 

Now we are in a position to prove Theorem 1, and we do it by noting it
is a part of the following.

\proclaim Theorem 11. Let $K$ have finite $u$ invariant. Let $K'/K$ be a any field 
extension of degree dividing $n$. 
For all degrees $d$, there is a symbol bound for $H^d(K',\mu_2)$ that only depends on 
$n$. Also, there is a limitation bound for $K'$ that only depends on $n$. 
Moreover, given a presentation of zero $\beta \in H^d(G,\mu_2)$ and $G = \Gal(L'/K')$, 
this limitation bound is realized by $d-1$ abelian $G$ Galois groups. 

\proof As we have said all along, we prove this by induction and so we assume this statement 
for all degrees $j < d$. By Theorem 5 we have the symbol boundedness of $K'$ in degree $d$. 
By Corollary 10, we have the limitation bound in degree $d$.~\qed 

\bigskip

\leftline{\bf{References}}
\medskip
[E] Evens, Leonard ``The Cohomology of Groups'', Clarendon Press 
Oxford/New York/Tokyo 1991 
\medskip
[NSW] Neukirch, Jurgen; Schmidt, Alexander; Wingber, Kay; 
``Cohomology of Number Fields'' Springer Berlin Heidelberg 2000 
\medskip
[OVV] Orlov, D.; Vishik, A,; Voevodsky, V. {\it An Exact Sequence for 
$K_*^M/2$ with Applications to Quadratic Forms}, Annals of Math. {\bf 165} 
(2007), 1-13 
\medskip
[Sa] Saltman, David J. {\it Generic Galois Extensions and Problems in Field Theory}, 
Adv in Math. {\bf 43}, 250-283 (1982)
\medskip
[Sc] Sharlau, Winfried ``Quadratic and Hermitian Forms'', Springer-Verlag, 
Berline/Heidelberg/New York/Tokyo 1985
\medskip
[Sw] Swan, R.G. {\it Invariant Rational Functions and a Problem of Steenrod}, 
Invent. Math. {\bf 7} (1969) p. 148-158
\end

which has the form 
$\Ind_{A_1}^A(K_1)$ where $K_1/F_1$ is $A_1$ Galois and $K_1$ is a field. 
Then 
The $A$ map $S \to L_1$ we also call $\phi$. 
Let $\beta_1 \in H^i(A_1,\mu_2)$ be the restriction of $\beta$. 
Then $K_1/F_1$ defines an inflation map $H^i(A_1,\mu_2) \to H^i(F_1,\mu_2)$ 
and $\beta_1$ has an image, $\alpha_1$, we call a {\bf specialization} 
of $\alpha \in H^i(F,\mu_2)$. 
Moreover we get a presentation of $\alpha_1$ as follows. 
There is a canonical $A_1$ module projection $L_1 \to K_1$ 
under which $L^{A_1}$ maps to $F_1$. We can choose $R_1$ with 
$q(R_1) = F_1$ and $R_1$ containing the image of $S^{A_1}$ 
(so $R_1$ also ccontains $\phi(R)$). Having done this, we set 
$S_1 = S \otimes_{\phi} R_1$. Then $S_1/R_1$ is $A_1$ Galois 
and together with $A_1$ and $\beta_1$ define a presentation of 
$\alpha_1$ which is a specialization 
of the presentation of $\alpha$. After we talk about dense specialization 
below it will be clear the choice of $S/R$ (inside a fixed $L/F$) is unimportant and so we will describe a presentation 
as being determined by $\beta \in H^i(A,\mu_2)$ and the $A$ Galois $L/F$.